\documentclass[12pt]{amsart}
\usepackage[utf8]{inputenc}
\usepackage{array}
\usepackage{subfiles}
\usepackage{hyperref}
\usepackage{cleveref}
\usepackage{amsmath}
\usepackage{amssymb}
\usepackage{amsthm}
\usepackage{mathtools}
\usepackage{thmtools}
\usepackage{thm-restate}
\usepackage[margin=1in]{geometry}
\usepackage{cleveref}
\usepackage{ytableau}
\usepackage{comment}
\usepackage{soul}
\usepackage{tikz}
\usetikzlibrary{positioning}
\usepackage{tabularx}
\usepackage{graphicx}
\usepackage{graphics}
\usepackage{caption}
\usepackage{float}

\usepackage[T1]{fontenc} 

\usepackage[english]{babel}

\numberwithin{equation}{section}

\theoremstyle{definition}
\newtheorem{theorem}{Theorem}[section]
\newtheorem*{theorem*}{Theorem}
\newtheorem{example}[theorem]{Example}
\newtheorem*{example*}{Example}
\newtheorem{lemma}[theorem]{Lemma}
\newtheorem*{lemma*}{Lemma}
\newtheorem{corollary}[theorem]{Corollary}
\newtheorem*{corollary*}{Corollary}
\newtheorem{definition}[theorem]{Definition}
\newtheorem*{definition*}{Definition}

\newtheorem*{proposition*}{Proposition}

\newtheorem*{remark*}{Remark}
\newtheorem{conjecture}[theorem]{Conjecture}

\usepackage{ mathrsfs }

\title{Fertility Numbers of Consecutive $S_3$ Pattern-Avoiding Stack-Sorting maps}

\author{Jurgis Kemeklis}\address{\textsc{J. Kemeklis}, Phillips Exeter Academy, NH, 03833} \email{jurgis.kemeklis@gmail.com}

\begin{document}

\maketitle

\begin{abstract}
In this paper, we show that for all length 3 patterns, all positive integers are fertility numbers for the consecutive-pattern-avoiding stack-sorting map $\textrm{SC}_\sigma$, which resolves a conjecture by Defant and Zheng \cite{defant}. The paper ends with a conjecture. 
\end{abstract}

\section{Introduction} 
\label{intro}

\textbf{Background:} In 1968, Knuth \cite{knuth} started investigating pattern avoidance in permutations by defining a stack-sorting machine. He ended up showing that for some input $\tau \in S_n$ of this machine, the output is $12\cdots n$ if and only if the permutation $\tau$ doesn't have the pattern $231$. A deterministic version of this machine was introduced in 1990 by West \cite{west}. He defined it to be a map $s$ on permutations of length $n$, and it was an identical version of Knuth's machine except it had an additional restriction: it avoided the relative order $21$ in the stack. Analogously to Knuth's machine, $s$ sorted a permutation to the identity if and only if the permutation avoided the pattern $231$. Not too long ago, Albert, Homberger, Pantone, Shar, and Vatter \cite{albert} introduced a modified version of West's machine, the $C$ machine, in which the stack had the same relative order as some element of the permutation class $C$ when reading it from the top to the bottom. Cerbai, Claesson, and Ferrari recently provided a generalization of $s$, namely the machine $s_\sigma$, which, by working in a right-greedy manner, avoids the case for the numbers in the stack to have relative order $\sigma$ when read from the top to the bottom. By this definition of $s_\sigma$, West explored the map $s_{21}$. Although the article \cite{cerbai_restrictedstacks} is recent, it has already inspired several subsequent papers \cite{baril, berlow, cerbai_2021_sorting, cerbai_claesson_anders_ferrari}.

\begin{figure}[h]\begin{center}
\begin{tikzpicture}[scale=0.5]
\draw[thick] (0,0) -- (2,0) -- (2,-2.5) -- (3,-2.5) -- (3,0) -- (5,0);
\node[fill = white, draw = white, scale = 0.8] at (4,.5) {52413};
\node[fill = white, draw = white, scale = 0.8] at (2.5,-2) {};
\node[fill = white, draw = white, scale = 0.8] at (2.5,-1.3) {};
\node[fill = white, draw = white, scale = 0.8] at (6,-1) {$\rightarrow$};
\end{tikzpicture}
 \begin{tikzpicture}[scale=0.5]
\draw[thick] (0,0) -- (2,0) -- (2,-2.5) -- (3,-2.5) -- (3,0) -- (5,0);
\node[fill = white, draw = white, scale = 0.8] at (4,.5){2431};
\node[fill = white, draw = white, scale = .8] at (2.5,-2) {5};
\node[fill = white, draw = white, scale = 0.8] at (.5,.5) {};
\node[fill = white, draw = white, scale = 0.8] at (6,-1) {$\rightarrow$}; \end{tikzpicture}
\begin{tikzpicture}[scale=0.5]
\draw[thick] (0,0) -- (2,0) -- (2,-2.5) -- (3,-2.5) -- (3,0) -- (5,0);
\node[fill = white, draw = white, scale = 0.8] at (4,.5){431};
\node[fill = white, draw = white, scale = .8] at (2.5,-2) {5};
\node[fill = white, draw = white, scale = 0.8] at (2.5, -1.3){2};
\node[fill = white, draw = white, scale = 0.8] at (6,-1) {$\rightarrow$};
\end{tikzpicture}
\begin{tikzpicture}[scale=0.5]
\draw[thick] (0,0) -- (2,0) -- (2,-2.5) -- (3,-2.5) -- (3,0) -- (5,0);
\node[fill = white, draw = white, scale = 0.8] at (4,.5){413};
\node[fill = white, draw = white, scale = 0.8] at (2.5,-2) {5};
\node[fill = white, draw = white, scale = 0.8] at (.5,.5) {2};
\node[fill = white, draw = white, scale = 0.8] at (6,-1) {$\rightarrow$}; \end{tikzpicture}
\begin{tikzpicture}[scale=0.5]
\draw[thick] (0,0) -- (2,0) -- (2,-2.5) -- (3,-2.5) -- (3,0) -- (5,0);
\node[fill = white, draw = white, scale = 0.8] at (4,.5){13};
\node[fill = white, draw = white, scale = 0.8] at (2.5,-2) {5};
\node[fill = white, draw = white, scale = 0.8] at (2.5,-1.3) {4};
\node[fill = white, draw = white, scale = 0.8] at (.5, .5) {2};
\node[fill = white, draw = white, scale = 0.8] at (6,-1) {$\rightarrow$}; \end{tikzpicture}\begin{tikzpicture}[scale=0.5]
\draw[thick] (0,0) -- (2,0) -- (2,-2.5) -- (3,-2.5) -- (3,0) -- (5,0);
\node[fill = white, draw = white, scale = 0.8] at (4,.5){3};
\node[fill = white, draw = white, scale = 0.8] at (2.5,-2) {5};
\node[fill = white, draw = white, scale = 0.8] at (2.5,-1.3) {4};
\node[fill = white, draw = white, scale = 0.8] at (2.5,-.6) {1};
\node[fill = white, draw = white, scale = 0.8] at (.5,.5) {2};
\node[fill = white, draw = white, scale = 0.8] at (6,-1) {$\rightarrow$}; \end{tikzpicture}\begin{tikzpicture}[scale=0.5]
\draw[thick] (0,0) -- (2,0) -- (2,-2.5) -- (3,-2.5) -- (3,0) -- (5,0);
\node[fill = white, draw = white, scale = 0.8] at (4,.5){3};
\node[fill = white, draw = white, scale = 0.8] at (2.5,-2) {5};
\node[fill = white, draw = white, scale = 0.8] at (2.5,-1.3) {4};
\node[fill = white, draw = white, scale = 0.8] at (.5,.5) {21};
 \node[fill = white, draw = white, scale = 0.8] at (6,-1) {$\rightarrow$}; \end{tikzpicture}\begin{tikzpicture}[scale=0.5]
\draw[thick] (0,0) -- (2,0) -- (2,-2.5) -- (3,-2.5) -- (3,0) -- (5,0);
\node[fill = white, draw = white, scale = 0.8] at (4,.5){};
\node[fill = white, draw = white, scale = 0.8] at (2.5,-2) {5};
\node[fill = white, draw = white, scale = 0.8] at (2.5,-1.3) {4};
\node[fill = white, draw = white, scale = 0.8] at (2.5,-.6) {3};
\node[fill = white, draw = white, scale = 0.8] at (.5,.5) {21};
 \node[fill = white, draw = white, scale = 0.8] at (6,-1) {$\rightarrow$}; \end{tikzpicture}
 \begin{tikzpicture}[scale=0.5]
\draw[thick] (0,0) -- (2,0) -- (2,-2.5) -- (3,-2.5) -- (3,0) -- (5,0);
\node[fill = white, draw = white, scale = 0.8] at (2.5,-1.3) {4};
 \node[fill = white, draw = white, scale = 0.8] at (2.5,-2) {5};
\node[fill = white, draw = white, scale = 0.8] at (.5,.5) {213};
\node[fill = white, draw = white, scale = 0.8] at (6,-1) {$\rightarrow$};
 \end{tikzpicture}
\begin{tikzpicture}[scale=0.5]
\draw[thick] (0,0) -- (2,0) -- (2,-2.5) -- (3,-2.5) -- (3,0) -- (5,0);
 \node[fill = white, draw = white, scale = 0.8] at (2.5,-2) {5};
\node[fill = white, draw = white, scale = 0.8] at (.5,.5) {2134};
\node[fill = white, draw = white, scale = 0.8] at (6,-1) {$\rightarrow$};
 \end{tikzpicture}
\begin{tikzpicture}[scale=0.5]
\draw[thick] (0,0) -- (2,0) -- (2,-2.5) -- (3,-2.5) -- (3,0) -- (5,0);
\node[fill = white, draw = white, scale = 0.8] at (.5,.5) {21345};\end{tikzpicture}
\end{center}
\caption{The stack-sorting map $SC_{213}$ on $\tau = 52413$.}
\label{fig1}
\end{figure}
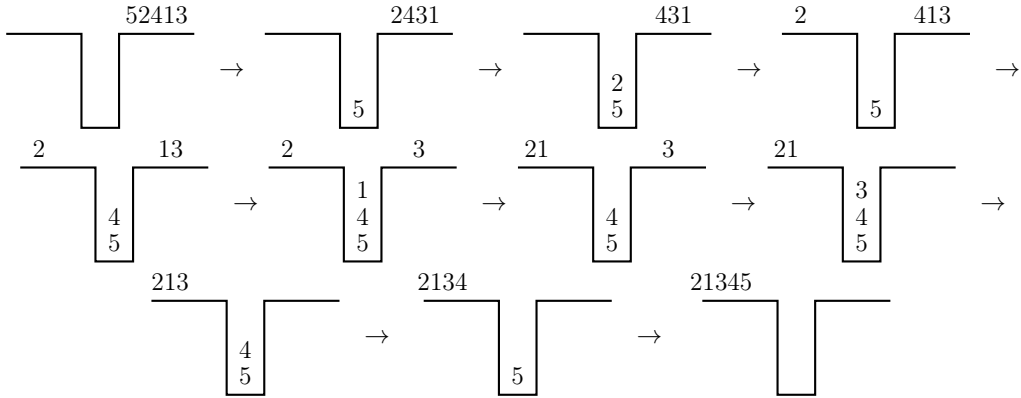

Defant and Zheng’s paper \cite{defant} introduced consecutive-pattern-avoiding stack-sorting maps \( \textrm{SC}_\sigma \). Given an input permutation \( \tau = \tau_1 \cdots \tau_n \), the procedure considers the permutation obtained by reading the stack's contents from top to bottom. If placing (but not yet \textit{pushing}) the next input permutation entry on top of the stack results in a stack whose top three numbers have consecutive relative order \( \sigma \), then the entry at the top of the stack (the last \textit{pushed} number) is moved to the end of the growing output permutation (this operation is called a \textit{pop}). In all other cases, the entry is \textit{pushed} onto the top of the stack. The procedure stops when the output permutation reaches length \( n \). For example, \Cref{fig1} shows that \( \textrm{SC}_{213}(52413) = 21345 \).

Defant and Zheng \cite{defant} raise a series of conjectures, two of which were resolved in a paper by Choi \cite{choi} and a paper by Seidel and Sun \cite{seidel}. In this paper, we resolve the conjecture 8.3, which states that for every $\sigma \in S_3$ and every positive integer $f$, there exists a permutation $\pi$ such that $|\textrm{SC}^{-1}_{\sigma}(\pi)| = f$. Traditionally, a nonnegative integer is called a \textit{fertility} number if it is equal to the number of preimages of a permutation under West’s stack-sorting map $s$. In this paper, we apply this concept of fertility to the consecutive-pattern-avoiding stack-sorting maps $\textrm{SC}_\sigma$. So the fertility of a permutation $\tau$ with $\sigma \in S_3$ is $|\textrm{SC}_\sigma^{-1}(\pi)|$. Reframing the conjecture 8.3 \cite{defant}, we find that we need to prove that:
\begin{theorem} \label{first+}
For every positive integer f and every $\sigma \in S_3$, the number f is a fertility number for $SC_\sigma$. \\
\end{theorem}

The rest of the paper is organized as follows. In \Cref{prelim}, we establish preliminaries. In \Cref{proofs}, we present the proofs of our main results. In \Cref{futuredirections}, we suggest potential future directions of study. 

\section{Preliminaries}
\label{prelim}
\subsection{Definitions} \label{definitions}

The following definitions come from Defant's and Zheng's paper \cite{defant}:\\
\begin{definition}
$S_n$ is the set of all permutations of $(1, 2, \cdots, n)$.
\end{definition}
Recall that in this article, a permutation is an ordering of the elements of the set $[n]$ for some $n$.
\begin{definition} If $\pi = \pi_1\pi_2\cdots\pi_n$ is a sequence of $n$ distinct integers, then the \textit{standardization} of $\pi$ is the permutation in $S_n$ obtained by replacing the $i$th-smallest entry in $\pi$ with $i$ for all $i$. \end{definition}
\begin{example} The \textit{standardization} of 4829 is 2314. \end{example} We say two sequences have the same relative
order if their \textit{standardizations} are equal. We say a permutation $\pi$ contains a permutation $\sigma$ as a pattern if there is a (not necessarily consecutive) subsequence of $\pi$ that has the same relative order as $\sigma$; otherwise, $\pi$ avoids $\sigma$. We say $\pi$ contains $\sigma$ consecutively if $\pi$ has a consecutive subsequence with the same relative order as $\sigma$; otherwise, $\pi$ avoids the consecutive pattern $\sigma$.
\begin{definition}
If $\pi \in S_n$, then $\textrm{rev}(\pi) = \pi_n \cdots \pi_1$ denotes the reverse of $\pi$.
\end{definition}
\begin{example}
$\textrm{rev}(3412) = 2143$.
\end{example}
\begin{definition}
If $\pi \in S_n$, then $\textrm{comp}(\pi) = (n+1 - \pi_1)(n+1 - \pi_{2})\cdots(n+1 - \pi_n)$ is the complement of $\pi$.
\end{definition}
\begin{example}
$\textrm{comp}(123) = 321$.
\end{example}

\begin{definition}
For $\pi = \pi_1\pi_2\cdots\pi_n \in S_n$ $\textrm{ind}_\pi(x)$ is a function that returns the index of $x$ in $\pi$. We assume that the index starts at 1, and we also use the term position interchangeably.
\end{definition}
\begin{example}
If $\pi = 51243$, then $\textrm{ind}_\pi(2) = 3$.
\end{example}
Now let us define some other terms, not mentioned in previous papers.
\begin{definition}
$\textrm{CRO}(\textrm{SC}_\sigma(\pi))$ is the number of pops that happen due to Consecutive Relative Order being $\sigma$ while performing the $\textrm{SC}_\sigma(\pi)$ function.
\end{definition}
\begin{example}
$\textrm{CRO}(\textrm{SC}_{213}(52413)) = 2$ because in \Cref{fig1} you can see that the numbers 2 and 1 are popped out because of relative order being $\sigma$.
\end{example}
\begin{definition}
A \textit{combination} of the stack and the numbers left in the input $\tau$ is the sequence we get if we were to push the remaining numbers from the input $\tau$ to the stack and read it from the bottom to the top.
\end{definition}
\begin{example}
The combination in \Cref{fig2} is 543.
\end{example}
\begin{figure}[h]\begin{center}
\begin{tikzpicture}[scale=0.5]
\draw[thick] (0,0) -- (2,0) -- (2,-2.5) -- (3,-2.5) -- (3,0) -- (5,0);
\node[fill = white, draw = white, scale = 0.8] at (4,.5){3};
\node[fill = white, draw = white, scale = 0.8] at (2.5,-2) {5};
\node[fill = white, draw = white, scale = 0.8] at (2.5,-1.3) {4};
\node[fill = white, draw = white, scale = 0.8] at (.5,.5) {21};
\end{tikzpicture}
\end{center}
\caption{An instance of $\textrm{SC}_{213}(52413)$.}
\label{fig2}
\end{figure}
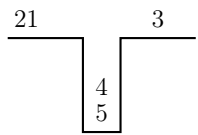
\begin{definition}
The \textit{fertility} of a permutation $\pi$ with $\sigma \in S_3$ is $|\textrm{SC}_\sigma^{-1}(\pi)|$.
\end{definition}
\begin{example}
The fertility of $\pi = 1243$ with $\sigma = 213$ is 2 because $|\textrm{SC}_{213}^{-1}(1243)| = 2$.
\end{example}
Whenever we write $\sigma$ in the rest of the paper assume that $\sigma \in S_3$. Whenever we write $\tau$ and $\pi$, assume that they are the input and output of $\textrm{SC}_\sigma$, respectively.

\section{Proofs of the Main Results}
\label{proofs}
We present a chart below illustrating a summarized version of our main findings.
\begin{center}
\begin{tabular}{ |c|l|c|c|}

\hline
$\sigma$ & Permutation $\pi$ & Fertility number& Theorem \\
\hline
$123$ & $n(n-1)(n-2) \dots321(n+1)$ & $n$ & \Cref{theorem3}\\
\hline
$321$ & $234\cdots(n-1)n(n+1)1$& $n$ &\Cref{theorem3}\\
\hline
$312$ & $123\cdots (n-1)(n+1)n$& $n$ &\Cref{theorem4} \\
\hline
$132$ & $(n+1)n(n-1)\cdots312$& $n$ & \Cref{theorem4}\\
\hline
$213$ & $123\cdots (n-4)(n-2)(n-3)(n-1)n(n+1)$& $n-1$ if $n\geq6$ &\Cref{theorem5} \\
\hline
$231$ & $(n+1)n(n-1)\cdots 645321$& $n-1$ if $ n\geq6$ & \Cref{theorem5}\\

\hline

\end{tabular}
\end{center}

\begin{lemma}\label{lemma1}
If $\textrm{SC}_\sigma(\tau) = \pi_1\pi_2\cdots\pi_n$, then the number $\pi_i$ is popped during the $i^{\textrm{th}}$ pop of the stack.
\end{lemma}
\begin{proof}
By the definition of the $\textrm{\textrm{SC}}_\sigma(\tau) = \pi$ function, the only way elements of the stack go to the output is once they are popped out. Since the elements of the output are $\pi = \pi_1\pi_2, \cdots \pi_n$, the $i^{\textrm{th}}$ pop has to be the number $\pi_i$.
\end{proof}

\begin{lemma}\label{lemma5}
Suppose that $\tau \in S_{n+1}$ and $ \sigma \in S_3$. If $\textrm{\textrm{SC}}_\sigma(\tau) = \pi_1\pi_2\cdots\pi_{n+1}$, then we get that $ \tau_1 = \pi_{n+1}$.
\end{lemma}
\begin{proof}
By \Cref{lemma1}, the $(n+1)^{\textrm{th}}$ number of the input will be popped on the $(n+1)^{\textrm{th}}$ pop, which in this case is the very last pop of the stack. From this, we can deduce that $\tau_1 = \pi_{n+1}$.
\end{proof}
Now we will prove a lemma concerning the value of $\textrm{\textrm{CRO}}({\textrm{SC}}_\sigma(\tau))$.
\begin{lemma} \label{lemma2}
After $\textrm{\textrm{CRO}}(\textrm{\textrm{SC}}_\sigma(\tau))$ pops of the stack, the combination of the stack and $\tau$ will be popped in reverse order to the output.
\end{lemma}

\begin{proof}
Notice that after the last pop, due to the relative order being $\sigma$, we will have the rest of the numbers of $\tau$ being pushed to the stack, only to be popped out from the stack to the output in reverse order.
\end{proof}
\begin{lemma} \label{lemma3}
If $\textrm{SC}_\sigma(\tau) = \pi$, then $\textrm{CRO}(\textrm{SC}_\sigma(\textrm{rev}(\pi)) = 0$ if and only if $\textrm{SC}_\sigma(\textrm{rev}(\pi)) = \pi$.
\end{lemma}
\begin{proof}

Let us prove the first direction of this lemma. Having $\textrm{CRO}(\textrm{SC}_\sigma(\textrm{rev}(\pi)) = 0$ means that the output of $\textrm{SC}_\sigma(\textrm{rev}(\pi))$ is just the reverse of the input $\textrm{rev}(\pi)$ by \Cref{lemma2}, since all the numbers are just pushed into the stack, and once there are no more numbers to push, the numbers are popped out one by one from the top. Notice that $\textrm{rev}(\textrm{rev}(\pi)) =\pi$, which gives us our desired result.
For the second direction, let us assume otherwise. For the sake of contradiction, suppose that $\textrm{SC}_\sigma(\textrm{rev}(\pi)) = \pi$ and $ \textrm{CRO}(\textrm{SC}_\sigma(\textrm{rev}(\pi)) > 0$ but then let us fix a moment in time, while the numbers are being pushed into the stack and then some number is popped first due to relative order being $\sigma$. Notice that then the first number in the output $\pi$ is just the last number in $\textrm{rev}(\pi)$, which at this current moment of time has not even been pushed into the stack. We reach a contradiction, as sought.

\end{proof}

\begin{theorem} \label{theorem1}
Suppose we have $\sigma \in S_3$ and $\pi, \tau \in S_n$. If $SC_\sigma(\tau) = \pi$, then \\$SC_{\textrm{comp}(\sigma)}(\textrm{comp}(\tau)) = \textrm{comp}(\pi)$.
\end{theorem}
\begin{proof}
Basic symmetry applies here. At any time that relative order $\sigma$ appears in the stack while performing the $\textrm{SC}_\sigma(\tau)$ operations, is also the same point when $\textrm{comp}(\sigma)$ would appear when performing the $\textrm{SC}_{\textrm{comp}(\sigma)}(\textrm{comp}(\tau))$ operations. So the way the numbers are popped is identical to $SC_\sigma(\tau)$, except the numbers are not $\pi = \pi_1\pi_2\cdots\pi_n$ but they are $\textrm{comp}(\pi) = (n+1-\pi_1)(n+1-\pi_2)\cdots(n+1-\pi_n)$.
\end{proof}
\begin{corollary}\label{cor1}
    Suppose we have $\sigma \in S_3$ and $\pi \in S_n$, such that $|\textrm{SC}^{-1}_{\sigma}(\pi)| = f$, for some positive integer $f$. We then also have that $| \textrm{SC}^{-1}_{\textrm{comp}(\sigma)}(\textrm{comp}(\pi))| = f$.
\end{corollary}
\begin{proof}
    It follows directly from \Cref{theorem1}.
\end{proof}
Using \Cref{cor1} it is enough to find constructions of $\pi$ for $\sigma = \{123, 213, 312\}$, because for $\sigma = \{321, 231, 132 \}$ we can just take the $\textrm{comp}(\pi)$ as a construction for the complement of $\{123, 213, 312\}$, respectively. Now let us denote some input $\tau \in S_{n+1}$ as $\tau = \tau_1\tau_2\dots \tau_{n+1}$. We will prove that all integers are fertility numbers for the maps $\textrm{\textrm{SC}}_{123}$ and $\textrm{\textrm{SC}}_{321}$.
\begin{theorem}\label{theorem3}

If $\sigma = 123$, then the fertility number of $\pi = n(n-1)(n-2) \dots21(n+1)$ is $n$.
\end{theorem}

\begin{proof}
Because $\pi_{n+1} = n+1$, we can apply \Cref{lemma5} to get $\tau_1 = n+1$.
We will prove this theorem by analyzing different cases for $\textrm{\textrm{CRO}}(\textrm{\textrm{SC}}_\sigma(\tau))$. More specifically, we will prove the following claim: if $\textrm{\textrm{CRO}}(\textrm{\textrm{SC}}_\sigma(\tau)) = k$ where $k \in \{0, 1, \ldots, n-1\}$, then there is exactly one such preimage $\tau$ for each specific value of $k$. Let us first start with the case when $\textrm{\textrm{CRO}}(\textrm{\textrm{SC}}_\sigma(\tau)) = 0$. By \Cref{lemma2}, $\tau$ is just the reverse order of the output $\pi$, and so we find that there is exactly one preimage in this case, namely: \begin{align*}\tau = (n+1)12\cdots(n-2)(n-1)n.\end{align*}
Now assume that \begin{align*}\textrm{\textrm{CRO}}(\textrm{\textrm{SC}}_\sigma(\tau)) = m \in \{1, 2, \ldots, n-1\}.\end{align*}
Let us say $i$ pops have currently happened, where $0\leq i<m$. By \Cref{lemma1}, the numbers $n, n-1, \ldots, n-i+1$ have been popped out. Notice that $n-i$ is the second largest number (smaller only to $n+1$) in the combination of the stack and numbers left in input at this moment in time. Combining this observation with \Cref{lemma1} stating that $n-i$ will be popped on the $(i+1)^{\textrm{th}}$ pop due to relative order being $\sigma$ implies that $n-i$ will be adjacent to the number $n+1$ which is at the bottom of the stack. Otherwise, $n-i$ would not be popped due to the relative order being $\sigma$. Repeating this line of reasoning for all possible values of $i$ from $0$ to $m-1$, we get that the input $\tau$ must look like this:\begin{align*} \tau = (n+1)n\cdots(n-m+1)\tau_{m+2}\tau_{m+3} \cdots \tau_{n+1}. \end{align*}

After $m$ pops, the combination of the stack and the numbers left in the input will look like this:
\begin{align*}(n+1)\tau_{m+2}\tau_{m+3} \cdots \tau_{n+1}.\end{align*} 
By \Cref{lemma2}: \begin{align*} (n-m)(n-m-1)\cdots1(n+1) &= \pi_{m+1}\pi_{m+2}\cdots \pi_{n+1} \\ &= \textrm{rev}\left((n+1)\tau_{m+2}\tau_{m+3} \cdots \tau_{n+1}\right) \\ &= \tau_{n+1} \cdots \tau_{m+3}\tau_{m+2}(n+1).\end{align*}
From this, we get that: \begin{align*} \tau_{n+1} = n-m, \ldots, \tau_{m+3} = 2, \tau_{m+2} = 1.\end{align*}
And so we prove that there is exactly one preimage for each specific value of $k$. This preimage is: \begin{align*}
\tau = (n+1)n\cdots(n-m+1)12\cdots (n-m).
\end{align*}
Since $\textrm{\textrm{CRO}}(\textrm{\textrm{SC}}_\sigma(\tau))$ can obtain all values from 0 to $n-1$, we find that there are exactly $n$ total preimages:
\begin{align*}
(n+1)123 \cdots (n-2)(n-1)n\\
(n+1)n123\cdots (n-2)(n-1) \\
(n+1)n(n-1)123\cdots (n-2)\\
\multispan{2}\hfil\vdots\hfil \\
(n+1)n(n-1)(n-2)\cdots321.
\end{align*}
\end{proof}
We will use a slightly different proof to prove that all integers are fertility numbers for the maps $\textrm{\textrm{SC}}_{312}$ and $\textrm{\textrm{SC}}_{132}$.

\begin{theorem} \label{theorem4}
If $\sigma = 312$, then the fertility number of $\pi = 12\cdots(n-1)(n+1)n$ is $n$.
\end{theorem}
\begin{proof}

Because $\pi_{n+1} = n$, we can apply \Cref{lemma5} to get $\tau_1 = n$.
We will prove this theorem by analyzing different cases for $\textrm{\textrm{CRO}}(\textrm{\textrm{SC}}_\sigma(\tau))$. More specifically, we will prove the following claim: if $\textrm{\textrm{CRO}}(\textrm{\textrm{SC}}_\sigma(\tau)) = k$ where $k \in \{0, 1, \ldots, n-1 \}$, then there is exactly one such preimage $\tau$ for each specific value of $k$. Let us first start with the case when $\textrm{\textrm{CRO}}(\textrm{\textrm{SC}}_\sigma(\tau)) = 0$. By \Cref{lemma2}, $\tau$ is just the reverse order of the output $\pi$, and so we find that there is exactly one preimage in this case, namely: \begin{align*}\tau = n(n+1)(n-1)\cdots21.\end{align*}
Now assume that \begin{align*}\textrm{\textrm{CRO}}(\textrm{\textrm{SC}}_\sigma(\tau)) = m \in \{1, 2, \ldots, n-1\}.\end{align*}
After $m$ pops, the combination of the stack and the numbers left in the input will look like this:
\begin{align*}n\tau_{a_1}\tau_{a_2} \cdots \tau_{a_{n-m}},\end{align*} where ($\tau_{a_1}\tau_{a_2} \cdots \tau_{a_{n-m}})$ is just some permutation of $(m+1, m+2, \ldots, n-1, n+1)$.
By \Cref{lemma2}: \begin{align*} (m+1)(m+2)\cdots(n-1)(n+1)n &= \pi_{m+1}\pi_{m+2}\cdots \pi_{n+1} \\ &= \textrm{rev}\left(n\tau_{a_1}\tau_{a_2} \cdots \tau_{a_{n-m}}\right) \\ &= \tau_{a_{n-m}} \cdots \tau_{a_2}\tau_{a_1}n.\end{align*}
From this, we get that: \begin{align*} \tau_{a_{n-m}} = m+1, \tau_{a_{n-m-1}} = m+2, \ldots, \tau_{a_2} = n-1, \tau_{a_1} = n+1. \end{align*} For any $i \in \{1, 2, \ldots, m \}$, we will prove the following claim:

\begin{equation} \label{eq1} 2 \leq \textrm{ind}_\tau(i) \leq m+1. \end{equation} The first inequality is obvious since $\tau_1 = n$. For the sake of contradiction, assume that for some $i \in \{1, 2, \ldots, m\}$ we have

\begin{align*} \textrm{ind}_\tau(i) > m+1. \end{align*}
Let us take maximal such $i$. Notice then that if $\textrm{ind}_\tau(i) >\textrm{ind}_\tau(\tau_{a_{n-m}})$, then the number $i$ will not be popped due to the relative order being $\sigma$. It follows that there exist such $j \in \{1, 2, \ldots, n-m-1\}$ that \begin{align*} \textrm{ind}_\tau(\tau_{a_j}) < \textrm{ind}_\tau(i) < \textrm{ind}_\tau(\tau_{a_{j+1}}). \end{align*} By the maximality of $i$ and the fact that $\tau_{a_j} > \tau_{a_{j+1}}$, we get that the number $i$ will not get popped due to the relative order being $\sigma$. We arrive at a contradiction as desired. From \Cref{eq1}, it follows that:
\begin{align*} \tau = n\tau_2\tau_3\cdots\tau_{m+1}(n+1)(n-1)\cdots(m+1). \end{align*} Now we will prove another claim: \begin{equation}\label{eq2}\tau_1 > \tau_2>\tau_3>\cdots>\tau_{m+1}. \end{equation} The first inequality is obvious. Now notice if $\tau_2<\tau_3$, then the relative order of $n,\tau_2,\tau_3$ in the stack looking from top to bottom will be $213$. We get that $\tau_2$ will be popped later than $\tau_3$, which is a contradiction to \Cref{lemma1}. It follows that $\tau_1 > \tau_2 > \tau_3$. Now, for the sake of contradiction, suppose that there exists some minimal $i \in \{3, 4, \ldots, m+1\}$ such that \begin{align*} \tau_1 > \tau_2>\tau_3> \cdots>\tau_{i-1} \end{align*} but $\tau_{i-1} < \tau_i$. Let us call $j \in \{2, 3, \ldots, i-1\}$ such index that \begin{align*} \tau_2> \tau_3> \cdots>\tau_{j+1} > \tau_i> \tau_j > \tau_{j+1} > \cdots > \tau_{i-1}.\end{align*}  It follows that the numbers $\tau_{i-1}, \tau_{i-2}, \ldots, \tau_{j+1}$ will be popped due to relative order being $\sigma$. But then the number $\tau_j$ will be popped later than $\tau_i$, which is a contradiction to \Cref{lemma1}. From \Cref{eq2}, we get that: $\tau_2 = m, \tau_3 = m-1, \ldots, \tau_{m+1} = 1$. And so we prove that there is exactly one preimage for each specific value of $k$. This preimage is: \begin{align*}
\tau = nm\cdots 21 (n+1)(n-1)\cdots (m+1).
\end{align*}
Since $\textrm{\textrm{CRO}}(\textrm{\textrm{SC}}_\sigma(\tau))$ can obtain all values from 0 to $n-1$, we find that there are exactly $n$ total preimages:

\begin{align*}
n(n+1)(n-1)(n-2)\cdots321\\
n1(n+1)(n-1)(n-2)\cdots32\\
n21(n+1)(n-1)(n-2)\cdots3\\
\multispan{2}\hfil\vdots\hfil \\
n(n-2)\cdots321(n-1)(n+1).\\
\end{align*}

\end{proof}
Lastly, we will prove that all positive integers greater than 4 are fertility numbers for the maps $\textrm{\textrm{SC}}_{213}$ and $\textrm{\textrm{SC}}_{231}$ \Cref{theorem5} in a different way than in the two solutions above. To do so, we will first prove another theorem that will be vital to the proof of \Cref{theorem5}.

\begin{theorem} \label{lemma4}
Suppose $\tau \in S_n$ and $\sigma = 213$. If we have that $\textrm{SC}_\sigma(\tau) = \pi = 12\cdots(n-4)(n-2)(n-3)(n-1)n(n+1)$ and $\textrm{CRO}(\textrm{SC}_\sigma(\tau)) = k$, where $k$ is some positive integer, then the first $k$ numbers of the output $\pi$ have the same relative order as they have in the input $\tau$.
\end{theorem}
\begin{proof}
 From the statement, we have that there will be exactly $k$ numbers popped due to the relative order being $\sigma$, and, by \Cref{lemma1}, these numbers will be: $\pi_1, \pi_2, \ldots, \pi_k$. From the statement we know that $\pi_1 = 1, \pi_2 = 2, \ldots, \pi_k = k$. For the sake of contradiction, let us assume that the numbers $\pi_1, \pi_2, \ldots, \pi_k$ do not have the same relative order as they have in the input. Then there exists some minimal index $i \in \{ 2, \ldots, k \}$ such that \begin{equation*} \textrm{ind}_\tau(\pi_1) < \textrm{ind}_\tau(\pi_2) < \cdots < \textrm{ind}_\tau(\pi_{i-2}) < \textrm{ind}_\tau(\pi_i) <  \textrm{ind}_\tau(\pi_{i-1}). \end{equation*} By \Cref{lemma1} we have that the number $\pi_{i-1}$ will be popped on the $(i-1)^{\textrm{th}}$ pop and the number $\pi_i$ will be popped on the $i^{\textrm{th}}$ pop. Since $\textrm{ind}_\tau(\pi_i) < \textrm{ind}_\tau(\pi_{i-1})$, we get that exactly after the $(i-1)^{\textrm{th}}$ pop if the next entry in the stack does not cause $\pi_i$ to be popped, then this number will not be popped on the $i^\textrm{th}$ pop which is a contradiction. It follows that \begin{equation}\label{eq3} \textrm{ind}_\tau(\pi_i) + 1 = \textrm{ind}_\tau(\pi_{i-1}).\end{equation} But after $i-2$ pops we get that the number $i-1$ is the smallest, and $i$ is the second smallest number in the combination. Combining \Cref{eq3} with this  means that the number $i-1$ will not be popped due to the relative order being $\sigma$ which produces contradiction as sought.
\end{proof}
Using \Cref{theorem1} on \Cref{lemma4}, we get that \Cref{lemma4} holds for $\sigma = 231$ as well. 
\begin{theorem} \label{theorem5}
If $\sigma = 213$, then the fertility number of $\pi = 12\cdots (n-4)(n-2)(n-3)(n-1)n(n+1)$ is $n-1$ for any positive integer $n \geq 6$.
\end{theorem}
\begin{proof}
From \Cref{lemma5}, we get that $\tau_1 = n+1$. Notice that the number $n-3$ is popped earlier than $n-2$ because of the relative order being $\sigma$ implies $\textrm{\textrm{SC}}_\sigma(\textrm{rev}(\pi)) \neq \pi$. By \Cref{lemma3} we get that \begin{align*} \textrm{\textrm{CRO}}(\textrm{\textrm{SC}}_\sigma(\textrm{rev}(\pi))) > 0 \implies \textrm{\textrm{CRO}}(\textrm{\textrm{SC}}_\sigma(\tau)) > 0 \end{align*} for any preimage of the output $\pi$. Let us denote the positions of the numbers $n$ and $n-1$ in the input $\tau$ by $i$ and $j$, respectively. In other words, $i = \textrm{ind}_\tau(n), j = \textrm{ind}_\tau(n-1)$. We will continue this proof by analyzing all of the possible cases for a permutation from $n+1$ numbers and showing that only $n-1$ of these permutations are the preimages to the output $\pi$. See \Cref{fig:TamB3} for a case diagram. Now let us analyze the first two cases:\\

\textit{Case 1:} if $j<i$, we can infer that all the numbers $(\tau_2, \tau_3 \cdots \tau_{j-1}, \tau_{j+1}, \tau_{j+2}, \ldots, \tau_{i-1})$ need to be popped due to relative order being $\sigma$. Otherwise, the number $n$ would be popped before the number $n-1$, which would be a contradiction. Another way of looking at this claim is by simply noting that all the numbers in between $\tau_1$ and $\tau_j$, as well as all the numbers in between $\tau_j$ and $\tau_i$, need to be popped due to their relative order being $\sigma$. Furthermore, notice that by \Cref{lemma1} the number $n-1$ needs to be popped on the $(n-1)^{\textrm{th}}$ pop, which means that \begin{align*} i = n+1 \implies \tau_{n+1} = n\end{align*} in this case. From the claim above and $\tau_{n+1} = n$, we can infer that all the numbers $\tau_2, \tau_3 \cdots \tau_{j-1}, \tau_{j+1}, \tau_{j+2}, \ldots, \tau_{n}$, which are just some permutation of 
\begin{align*}
\left( 1, 2, \ldots, n-4, n-3, n-2 \right)
\end{align*}
need to be popped due to relative order being $\sigma$. By \Cref{lemma4}, if there are exactly $n-2$ total pops due to the relative order being $\sigma$, we get that the first $n-2$ numbers of the output $\pi$ have the same relative order as in the input $\tau.$ Let us call $l$ the index of $n-3$, namely $l = \textrm{ind}_\tau(n-3)$. Let us prove that $l>j$: for the sake of contradiction, assume that $l<j$. We get that \begin{align*} \tau_{l-1} = n-2, \tau_{l-2} = n-4, \tau_{l-3} = n-5, \ldots, \tau_2 = 1 \implies l = n-1. \end{align*} But if \begin{align*} \tau = (n+1)12\cdots(n-5)(n-4)(n-2)(n-3)(n-1)n, \end{align*} then $ \textrm{\textrm{SC}}_\sigma(\tau) \neq \pi$. So we can assume that $l>j$, which just implies $ l=n$. Now let us prove that there are no preimages in this case by analyzing the position of $n-2$ in the permutation. Let us call $k = \textrm{ind}_\tau(n-2)$. If $k>j$, then $k = n-1$, but this is not possible since $n-3 < n-2$, and therefore $n-2$ would not be popped out due to the relative order being $\sigma$. If $k<j$ we get that \begin{align*}k = j-1 \implies \tau_{j-1} = n-2, \tau_{j-2} = n-4, \tau_{j-3} = n-5, \ldots, \tau_2 = 1 \end{align*} This situation can only happen when $j = n-1$, and so we get that \begin{align*}\tau = (n+1)123\cdots(n-4)(n-2)(n-1)(n-3)n \end{align*} but \begin{align*} \textrm{\textrm{SC}}_\sigma(\tau) \neq \pi\end{align*} which means that there are no preimages in this case as sought. Let us move on to the second case. \\
\textit{Case 2:} if $j>i$, then notice that since $n+1>n>n-1$ and $1 < i < j$ implies that none of the numbers $n+1, n, n-1$ will be popped due to the relative order being $\sigma$. Similarly to the first case, let us call $k = \textrm{ind}_\tau(n-2)$ and $l = \textrm{ind}_\tau(n-3)$. Now let us analyze two subcases for the number $n-3$ being popped out. First, let us suppose $n-3$ is popped out due to the relative order being $\sigma$. We will prove that there is exactly one preimage in this subcase. By \Cref{lemma1}, the number $n-3$ will be popped on the $(n-2)^{\textrm{th}}$ pop, which means that there will be $n-2$ total pops due to the relative order being $\sigma$. By \Cref{lemma4}, if there are exactly $n-2$ pops due to the relative order being $\sigma$, the first $n-2$ numbers of the output $\pi$ have the same relative order as in the input $\tau$, which in our case also means that $k<l$. Moreover, if $l<i$, then $\tau_{l-1} = n-2, \tau_{l-2} = n-4$. But notice that the numbers $\tau_l\tau_{l-1}\tau_{l-2}$ have a $123$ relative order, which means that the number $n-2$ would not be popped due to the relative order being $\sigma$. So henceforth, we can assume \begin{align*} l>i \implies l = n, j =n+1\implies \tau_n = n-3, \tau_{n+1} = n-1. \end{align*} Now if $k>i$, then by similar reasoning we get that the number $n-2$ would not get popped, so we can certainly claim that $k<i$, but by \Cref{lemma4} we get that \begin{align*}\tau_2 = 1, \tau_3 = 2, \ldots, \tau_{n-3} = n-4, \tau_{n-2} = n-2 \implies \tau_{n-1} = n.\end{align*} Putting \begin{align*} \tau =(n+1)123\cdots(n-4)(n-2)n(n-3)(n-1)\end{align*}  as the input does indeed give the desired result: $\textrm{\textrm{SC}}_\sigma(\tau) = \pi$, so there is exactly one preimage in this subcase as sought. Now let us assume that the number $n-3$ does not get popped due to the relative order being $\sigma$. We will prove there are exactly $n-2$ preimages in this subcase. By \Cref{lemma1} since the number $n-3$ is popped after the $(n-2)^{\textrm{th}}$ pop, we deduce that $l>j$. Notice that again, by \Cref{lemma1}, the number $n-2$ has to be popped on the $(n-3)^{\textrm{th}}$ pop. We claim that $k<j$. For the sake of contradiction, suppose that $k>j$. If $k>l$, we would get that the number $n-2$ is popped later than the number $n-3$ due to the relative order being $\sigma$, which is a contradiction to \Cref{lemma1}. If $k<l$, then once again the number $n-2$ is popped later than the number $n-3$, which is a contradiction to \Cref{lemma1}. So henceforth, we can assume $k<j$, but this implies that \begin{align*} l= n+1, j = n \implies \tau_{n+1} = n-3, \tau_n = n-1.\end{align*} Notice now that because of \Cref{lemma4}, the relative order of \begin{align*}(1,2, \ldots, n-4, n-2)\end{align*} is fixed. Combining this with the fact that we know the indexes of numbers $n+1, n-1$, and $n-3$, the number $n$ has $n-2$ possible positions. So, we get that there are a total of $1+(n-2)=n-1$ preimages:
\begin{align*}
(n+1)123\cdots(n-4)(n-2)\textcolor{blue}{n}\textcolor{red}{(n-3)}(n-1) \quad \quad & \text{1 preimage}\\
\left.
\begin{aligned}
& (n+1)123\cdots(n-4)(n-2)\textcolor{blue}{n}(n-1)\textcolor{red}{(n-3)}\\
& (n+1)123\cdots(n-4)\textcolor{blue}{n}(n-2)(n-1)\textcolor{red}{(n-3)}\\
&(n+1)123\cdots \textcolor{blue}{n}(n-4)(n-2)(n-1)\textcolor{red}{(n-3)}\\
& \multispan{2}\hfil\vdots\hfil \\
& (n+1)\textcolor{blue}{n}123\cdots(n-4)(n-2)(n-1)\textcolor{red}{(n-3)}
\end{aligned}
\right\} \quad & \text{$n-2$ preimages}
\end{align*}

\begin{center}
\begin{figure}[!hbtp]
\centering
\captionsetup{justification=centering} 
\scalebox{0.85}{
\begin{tikzpicture}
\node (ji) at (-4,10) {$j<i$};
\node (ij) at (4,10) {$j>i$};
\node (lj) at (-2, 8) {$l<j$};
\node (jl) at (-6, 8) {$l>j$};
\node (kj) at (-.5, 6) {$k<j$};
\node (jk) at (-3.5, 6) {$k>j$};
\node[color = red](empty) at (-6,6) {$\emptyset$};
\node[color = red] (empty2) at (-3.5,4) {$\emptyset$};
\node[color =red] (empty3) at (-0.5,4) {$\emptyset$};
\node (li) at (2, 8) {$l<i$};
\node[color = red] (empty4) at (2, 6) {$\emptyset$};
\node (il) at (6, 8) {$l>i$};
\node (ik2) at (4.5, 6) {$k>i$};
\node (ki2) at (7.5, 6) {$k<i$};
\node [color = green!70!black] (empty5) at (4.5, 4) {$1$};
\node (jl2) at (8.5, 4) {$l>j$};
\node (jk2) at (7.5, 2) {$k>j$};
\node (lk) at (6.5, 0) {$k>l$};
\node (kl) at (8.5, 0) {$k<l$};
\node (kj2) at (9.5, 2) {$k<j$};
\node[color = green!70!black](sol) at (10.5, 0) {$n-2$};
\node[color = red] (empty6) at (6.5, -2) {$\emptyset$};
\node[color = red] (empty7) at (8.5, -2) {$\emptyset$};

\node (start) at (0,12) {$i, j, k, l$};
\draw[] (start) -- (ji) -- (jl) -- (empty);
\draw (start) -- (ji) -- (lj) -- (jk) -- (empty2);
\draw (start) -- (ji) -- (lj) -- (kj) -- (empty3);
\draw (start) -- (ij) -- (li) -- (empty4);
\draw (start) -- (ij) -- (il) -- (ik2) --(empty5);
\draw (start) -- (ij) -- (il) -- (ik2) -- (il) -- (ki2) -- (jl2) -- (jk2) -- (lk) -- (empty6);
\draw (start) -- (ij) -- (il) -- (ik2) -- (il) --(ki2) -- (jl2) -- (jk2) -- (kl) -- (empty7);
\draw (start) -- (ij) -- (il) -- (ik2) -- (il) -- (ki2) -- (jl2) -- (kj2) -- (sol);

\end{tikzpicture}
}
\caption{A case diagram for \Cref{theorem5}}
\label{fig:TamB3}
\end{figure}
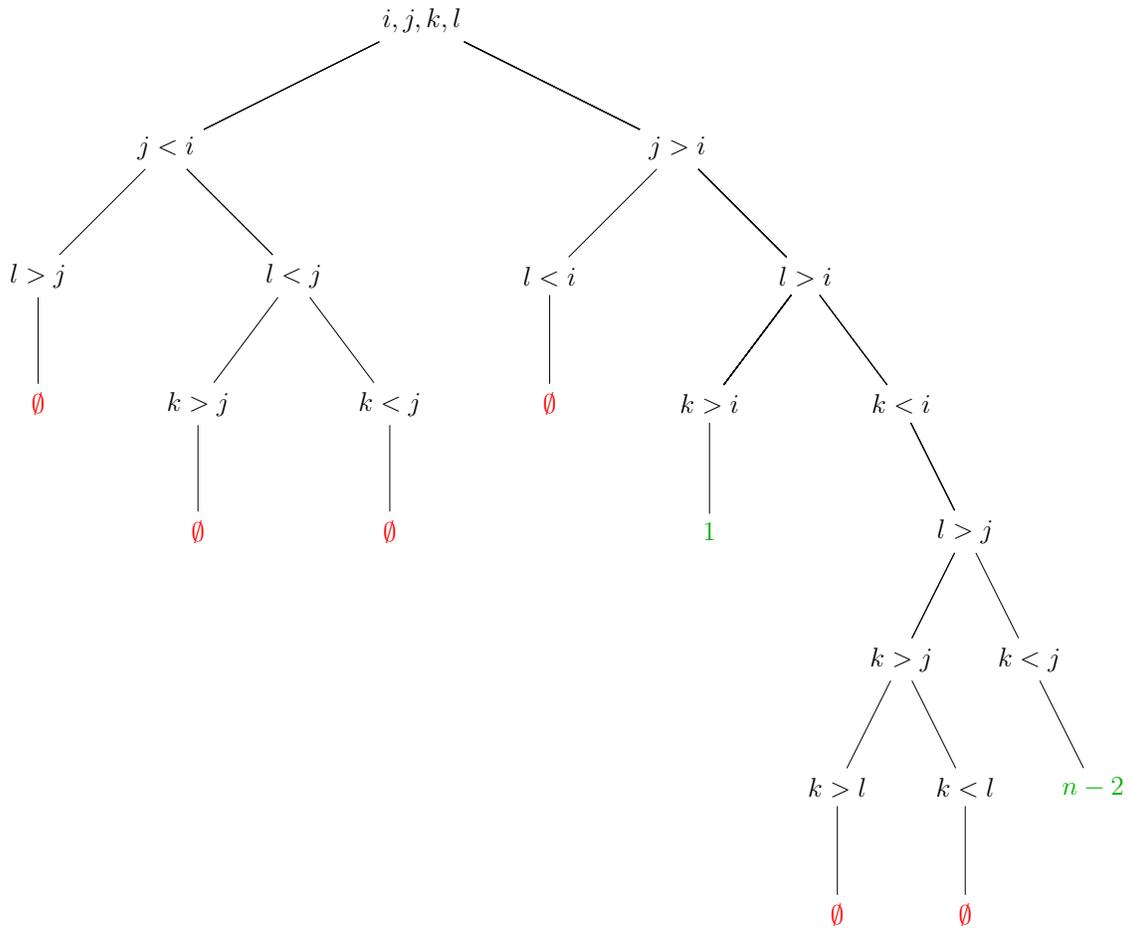
\end{center}

\end{proof}
Notice that \Cref{theorem5} gives all the fertility numbers $n \geq 5$ for $\sigma = 213$. Doing casework, we find that \begin{align*} |\textrm{\textrm{SC}}_{213}(4321)^{-1}| =1, |\textrm{\textrm{SC}}_{213}(1243)^{-1}| =2, |\textrm{\textrm{SC}}_{213}( 13524)^{-1}| =3, |\textrm{\textrm{SC}}_{213}(1234)^{-1}| =4,\end{align*} and so combining \Cref{theorem4}, \Cref{theorem3}, and \Cref{cor1} it follows that all positive integers are fertility numbers for the consecutive-pattern-avoiding stack-sorting map $\textrm{\textrm{SC}}_\sigma$ where $\sigma \in S_3$, proving \Cref{first+}. 


\section{Future Directions}
\label{futuredirections}

Analyzing all the permutations up to eight numbers with computer code, we raise the following conjecture.
\begin{conjecture} If $\sigma = \{213, 231\}$ and $A$ is the set of fertility numbers for all permutations out of $n$ numbers, then $\{1, 2, 3, \cdots, 2^{n-3}-1, 2^{n-3}\} \subset A$
    
\end{conjecture} 

\section*{Acknowledgements}
\noindent Thank you to Yunseo Choi, Alan Bu, and Mantas Bakšys for providing comments on the paper. 

\nocite{*}
\bibliographystyle{ieeetr}
\bibliography{partitions}

\begin{thebibliography}{10}

\bibitem{defant}
C.~Defant and K.~Zheng, ``Stack-sorting with consecutive-pattern-avoiding stacks,'' {\em Adv. in Appl. Math.}, vol.~128, pp.~102--192, 2021.

\bibitem{knuth}
D.~E. Knuth, {\em The Art of Computer Programming}.
\newblock Addison-Wesley Pub. Co, 1973.

\bibitem{west}
J.~West, {\em Permutations with Restricted Subsequences and stack-sortable Permutations}.
\newblock PhD thesis, 1990.

\bibitem{albert}
M.~H. Albert, C.~Homberger, J.~Pantone, N.~Shar, and V.~Vatter, ``Generating permutations with restricted containers,'' {\em J. Combin. Theory Ser. A}, vol.~157, pp.~205--232, 07 2018.

\bibitem{cerbai_restrictedstacks}
G.~Cerbai, A.~Claesson, and L.~Ferrari, ``Stack sorting with restricted stacks,'' {\em J. Combin. Theory Ser. A}, vol.~173, 07 2020.

\bibitem{baril}
J.-L. Baril, G.~Cerbai, C.~Khalil, and V.~Vajnovszki, ``Catalan and schröder permutations sortable by two restricted stacks,'' {\em Inform. Process. Lett.}, vol.~171, pp.~106138--106138, 10 2021.

\bibitem{berlow}
K.~Berlow, ``Restricted stacks as functions,'' {\em Discrete Math.}, vol.~344, 2021.

\bibitem{cerbai_2021_sorting}
G.~Cerbai, ``Sorting cayley permutations with pattern-avoiding machines *,'' {\em Australas. J. Combin.}, vol.~80, pp.~322--341, 2021.

\bibitem{cerbai_claesson_anders_ferrari}
G.~Cerbai, A.~Claesson, L.~Ferrari, and E.~Steingrímsson, ``Sorting with pattern-avoiding stacks: the 132-machine,'' {\em Electron. J. Combin.}, vol.~27, 08 2020.

\bibitem{choi}
Y.~Choi and Y.~Choi, ``Highly sorted permutations with respect to a 312-avoiding stack,'' {\em Enumer. Comb. Appl.}, vol.~3, 01 2022.

\bibitem{seidel}
I.~Seidel and N.~Sun, ``Periodic points of consecutive-pattern-avoiding stack-sorting maps,'' 08 2023.

\end{thebibliography}

\end{document}